\documentclass[11pt]{amsart}
\usepackage{stmaryrd}
\usepackage{mathrsfs}
\usepackage[centertags]{amsmath}
\usepackage{amsfonts, dsfont}
\usepackage{amssymb}
\usepackage{amsthm}
\usepackage[T2A,T1]{fontenc}
\usepackage[OT2,T1]{fontenc}
\DeclareSymbolFont{cyrletters}{OT2}{wncyr}{m}{n}
\DeclareMathSymbol{\Sha}{\mathalpha}{cyrletters}{"58}

\newcommand{\Zhe}{\mbox{\usefont{T2A}{\rmdefault}{m}{n}\CYRZH}}

\input xy
\xyoption{all}

\theoremstyle{plain}
\newtheorem{thm}{Theorem}[section]
\newtheorem{cor}[thm]{Corollary}
\newtheorem{lem}[thm]{Lemma}
\newtheorem{prop}[thm]{Proposition}

\theoremstyle{definition}
\newtheorem{defn}[thm]{Definition}
\newtheorem{remark}[thm]{Remark}
\newtheorem*{ack}{Acknowledgments}

\numberwithin{equation}{subsection}

\newcommand{\bd}{\begin{defn}}
\newcommand{\ed}{\end{defn}}
\newcommand{\bl}{\begin{lem}}
\newcommand{\el}{\end{lem}}
\newcommand{\bp}{\begin{prop}}
\newcommand{\ep}{\end{prop}}
\newcommand{\bt}{\begin{thm}}
\newcommand{\et}{\end{thm}}
\newcommand{\bc}{\begin{cor}}
\newcommand{\ec}{\end{cor}}
\newcommand{\br}{\begin{remark}}
\newcommand{\er}{\end{remark}}
\newcommand{\bdi}{\begin{diagram}}
\newcommand{\edi}{\end{diagram}}
\newcommand{\beq}{\begin{eqn}}
\newcommand{\eeq}{\end{eqn}}
\newcommand{\ba}{\begin{array}}
\newcommand{\ea}{\end{array}}
\newcommand{\bpf}{\begin{proof}}
\newcommand{\epf}{\end{proof}}

\newcommand{\Z}{\mathds{Z}}
\newcommand{\Q}{\mathds{Q}}
\newcommand{\Zp}{\mathds{Z}_{p}}
\newcommand{\Qp}{\mathds{Q}_{p}}
\newcommand{\al}{\alpha}
\newcommand{\Ga}{\Gamma}

\newcommand{\m}{\mathfrak{m}}

\DeclareMathOperator{\Sel}{Sel} \DeclareMathOperator{\Gal}{Gal}
\DeclareMathOperator{\Hom}{Hom} 
\DeclareMathOperator{\Cl}{Cl}

\newcommand{\ilim}{\displaystyle \mathop{\varinjlim}\limits}

\newcommand{\coker}{\mathrm{coker}\,}

\newcommand{\cyc}{\mathrm{cyc}}

\begin{document}

\title{The growth of fine Selmer groups}
\author{Meng Fai Lim}
\email{limmf@mail.ccnu.edu.cn}
\address{School of Mathematics and Statistics, Central China Normal University, No.152, Luoyu Road, Wuhan, Hubei 430079, CHINA}
\author{V. Kumar Murty}
\email{murty@math.toronto.edu}
\address{Department of Mathematics, University of Toronto, 40 St. George Street, Toronto, CANADA}
\date{\today}
\thanks{Research of VKM partially supported by an NSERC Discovery grant}
\maketitle

\begin{abstract}
 Let $A$ be an abelian variety defined over a number field
$F$. In this paper, we will investigate the growth of the $p$-rank
of the fine Selmer group in three situations. In particular, in
each of these situations, we show that there is a strong analogy
between the growth of the $p$-rank of the fine Selmer group and the
growth of the $p$-rank of the class groups.
\end{abstract}

\section{Introduction}

In the study of rational points on Abelian varieties, the Selmer group plays an important role. In Mazur's fundamental work \cite{Mazur}, the
Iwasawa theory of Selmer groups was introduced. Using this theory, Mazur was able to describe the growth of the
 size of the $p$-primary part of the Selmer group in $\Z_p$-towers. Recently, several authors have initiated the study of a certain subgroup, called the fine Selmer group. This subgroup, as well as the `fine' analogues of the Mordell-Weil group and Shafarevich-Tate group, seem to have stronger finiteness properties than the classical Selmer group (respectively,
Mordell-Weil or Shafarevich-Tate groups). The fundamental paper of Coates and Sujatha \cite{CS} explains some of these properties.
\medskip

Let $F$ be a number field and $p$ an odd prime. Let $A$ be an Abelian variety defined over $F$ and let $S$ be a finite set
of primes of $F$ including the infinite primes, the primes where $A$ has bad reduction, and the primes of $F$ over $p$.
Fix an algebraic closure $\overline{F}$ of $F$ and denote by $F_S$ the maximal subfield of $\overline{F}$ containing
$F$ which is unramified outside $S$. We set $G = \Gal(\overline{F}/F)$ and $G_S = \Gal(F_S/F)$.
\medskip

The usual $p^{\infty}$-Selmer group of $A$ is defined by
$$
 \Sel_{p^{\infty}}(A/F) = \ker\Big(H^1(G,A[p^{\infty}])\longrightarrow
\bigoplus_{v}H^1(F_v, A)[p^{\infty}]\Big).
$$
Here, $v$ runs through all the primes of $F$ and as usual, for a
$G$-module $M$, we write $H^*(F_v,M)$ for the Galois cohomology of
the decomposition group at $v$. Following \cite{CS}, the
$p^{\infty}$-fine Selmer group of $A$ is defined by
\[
R_{p^\infty}(A/F) = \ker\Big(H^1(G_S(F),A[p^{\infty}])\longrightarrow \bigoplus_{v\in S}H^1(F_v,
A[p^{\infty}])\Big).
 \]
 This definition is in fact independent of the choice of $S$ as can be seen from the exact sequence (Lemma \ref{indep of S})
 $$
 0 \longrightarrow R_{p^\infty}(A/F) \longrightarrow \Sel_{p^{\infty}}(A/F)
 \longrightarrow \bigoplus_{v|p}A(F_v)\otimes\Qp/\Zp.
 $$
\medskip

Coates and Sujatha study this group over a field $F_\infty$ contained in $F_S$ and for which $\Gal(F_\infty/F)$
is a $p$-adic Lie group. They set
$$
R_{p^\infty}(A/F_\infty)\ =\ \ilim R_{p^\infty}(A/L)
$$
where the inductive limit ranges over finite extensions $L$ of $F$
contained in $F_\infty$. When $F_\infty = F^{cyc}$ is the cyclotomic
$\Z_p$-extension of $F$, they conjecture that the Pontryagin dual
$Y_{p^\infty}(A/F_\infty)$ is a finitely generated $\Z_p$-module.
This is known {\em not} to be true for the dual of the classical
Selmer group. A concrete example of such is the elliptic curve
$E/\Q$ of conductor 11 which is given by
\[ y^2 +y=x^3 -x^2 -10 x- 20. \]
For the prime $p = 5$, it is known that the Pontryagin dual of the
$5^{\infty}$-Selmer group over $\Q^{\cyc}$ is not finitely generated
over $\Z_5$ (see \cite[\S 10, Example 2]{Mazur}). On the other hand
it is expected to be true if the Selmer group is replaced by a
module made out of class groups. Thus, in some sense, the fine
Selmer group seems to approximate the class group. One of the themes
of our paper is to give evidence of that by making the relationship
precise in three instances.
\medskip

Coates and Sujatha also study extensions for which $G_\infty = \Gal(F_\infty/F)$ is a $p$-adic Lie group of
dimension larger than $1$ containing $F^{cyc}$. They make a striking conjecture that the dual
$Y_{p^\infty}(A/F_\infty)$ of the fine Selmer group is pseudo-null as a module over the Iwasawa algebra
$\Lambda(G_\infty)$. While we have nothing to say about this conjecture, we investigate the growth of $p$-ranks of
fine Selmer groups in some pro-$p$ towers that are {\em not} $p$-adic analytic.
\medskip

\section{Outline of the Paper}

Throughout this paper, $p$ will always denote an odd prime. In the
first situation, we study the growth of the fine Selmer groups over
certain $\Zp$-extensions. It was first observed in \cite{CS} that
over a cyclotomic $\Zp$-extension, the growth of the fine Selmer
group of an abelian variety in a cyclotomic $\Zp$-extension exhibits
phenomena parallel to the growth of the $p$-part of the class groups
over a cyclotomic $\Zp$-extension. Subsequent papers \cite{A, JhS,
LimFine} in this direction has further confirmed this observation.
(Actually, in \cite{JhS, LimFine}, they have also considered the
variation of the fine Selmer group of a more general $p$-adic
representation. In this article, we will only be concerned with the
fine Selmer groups of abelian varieties.) In this paper, we will
show that the growth of the $p$-rank of fine Selmer group of an
abelian variety in a certain class of $\Zp$-extension is determined
by the growth of the $p$-rank of ideal class groups in the
$\Zp$-extension in question (see Theorem \ref{asymptotic compare})
and vice versa. We will also specialize our theorem to the
cyclotomic $\Zp$-extension to recover a theorem of Coates-Sujatha
\cite[Theorem 3.4]{CS}.
\medskip

In the second situation, we investigate the growth of the fine
Selmer groups over $\Z/p$-extensions of a fixed number field. We
note that it follows from an application of the Grunwald-Wang
theorem that the $p$-rank of the ideal class groups grows
unboundedly in $\Z/p$-extensions of a fixed number field. Recently,
many authors have made analogous studies in this direction replacing
the ideal class group by the classical Selmer group of an abelian
variety (see \cite{Ba, Br, Ce, Mat09}). In this article, we
investigate the analogous situation for the fine Selmer group of an
abelian variety, and we show that the $p$-rank of the fine Selmer
group of the abelian variety grows unboundedly in $\Z/p$-extensions
of a fixed number field (see Theorem \ref{class Z/p}). Note that the
fine Selmer group is a subgroup of the classical Selmer group, and
therefore, our results will also recover some of their results.
\medskip

In the last situation, we consider the growth of the fine Selmer
group in an infinite unramified pro-$p$ extensions. It is known that
the $p$-rank of the class groups is unbounded in such tower under
suitable assumptions. Our result will again show that we have the
same phenomenon for the $p$-rank of fine Selmer groups (see Theorem
\ref{Fine Sel in class tower}). As above, our result will also imply
some of the main results in \cite{LM, Ma, MO}, where analogous
studies in this direction have been made for the classical Selmer
group of an abelian variety.

\section{$p$-rank} \label{some cohomology lemmas}

In this section, we record some basic results on Galois cohomology
that will be used later. For an abelian group $N$, we define its
$p$-rank to be the $\Z/p$-dimension of $N[p]$ which we denote by
$r_p(N)$. If $G$ is a pro-$p$ group, we write $h_i(G) =
r_p\big(H^i(G,\Z/p)\big)$. We now state the following lemma which
gives an estimate of the $p$-rank of the first cohomology group.

\bl \label{cohomology rank inequalities} Let $G$ be a pro-$p$ group,
and let $M$ be a discrete $G$-module which is cofinitely generated
over $\Zp$.
 If $h_1(G)$ is finite, then $r_p\big(H^1(G,M)\big)$ is finite, and
 we have the following estimates for
 $r_p\big(H^1(G,M)\big)$
 \[
 h_1(G)r_p(M^G) -r_p \big( (M/M^G)^G\big)
\leq r_p\big(H^1(G,M)\big) \leq h_1(G)\big(\mathrm{corank}_{\Zp}(M)
+ \log_p(\big| M/M_{\mathrm{div}}\big|)\big).
\]
\el

\bpf
 See \cite[Lemma
 3.2]{LM}. \epf

We record another useful estimate.

\bl \label{estimate lemma}
 Let
 \[ W \longrightarrow X \longrightarrow Y \longrightarrow Z\]
 be an exact sequence of cofinitely generated abelian groups. Then we have
\[ \Big| r_p(X) - r_p(Y) \Big| \leq 2r_p(W) + r_p(Z).\]
\el

\bpf It suffices to show the lemma for the exact sequence
 \[ 0\longrightarrow W \longrightarrow X
 \longrightarrow Y \longrightarrow Z \longrightarrow 0.\]
We break up the exact sequence into two short exact sequences
 \[ 0\longrightarrow W \longrightarrow X
 \longrightarrow C \longrightarrow 0, \]
 \[ 0\longrightarrow C \longrightarrow Y
 \longrightarrow Z \longrightarrow 0.\]
 From these short exact sequences, we obtain two exact sequences of
finite dimensional $\Z/p$-vector spaces (since $W$, $X$, $Y$ and $Z$
are cofinitely generated abelian groups)
 \[ 0\longrightarrow W[p] \longrightarrow X[p]
 \longrightarrow C[p] \longrightarrow P \longrightarrow 0 , \]
 \[ 0\longrightarrow C[p] \longrightarrow Y[p]\longrightarrow Q
 \longrightarrow 0, \]
where $P\subseteq W/p$ and $Q\subseteq Z[p]$. It follows from these
two exact sequences and a straightforward calculation that we have
 \[ r_p(X) - r_p(Y) = r_p(W) - r_p(P) - r_p(Q). \]
 The inequality of the lemma is immediate from this.
 \epf
\medskip

\section{Fine Selmer groups} \label{Fine Selmer group section}

As before, $p$ will denote an odd prime. Let $A$ be an abelian
variety over a number field $F$. Let $S$ be a finite set of primes
of $F$ which contains the primes above $p$, the primes of bad
reduction of $A$ and the archimedean primes. Denote by $F_S$  the
maximal algebraic extension of $F$ unramified outside $S$. We will
write $G_S(F) = \Gal(F_S/F)$.
\medskip

 As stated in the introduction and following
\cite{CS}, the fine Selmer group of $A$ is defined by
\[ R(A/F) = \ker\Big(H^1(G_S(F),A[p^{\infty}])\longrightarrow \bigoplus_{v\in S}H^1(F_v,
A[p^{\infty}])\Big). \]
(Note that we have dropped the subscript $p^\infty$ on $R(A/F)$ as $p$ is fixed.)
\medskip

To facilitate further discussion, we also recall the definition of
the classical Selmer group of $A$ which is given by
\[ \Sel_{p^{\infty}}(A/F) = \ker\Big(H^1(F,A[p^{\infty}])\longrightarrow
\bigoplus_{v}H^1(F_v, A)[p^{\infty}]\Big), \] where $v$ runs through
all the primes of $F$. (Note the difference of the position of the
``$[p^{\infty}]$'' in the local cohomology groups in the
definitions.)
\medskip

At first viewing, it will seem that the definition of the fine
Selmer group depends on the choice of the set $S$. We shall show
that this is not the case.

\bl \label{indep of S}
 We have an exact sequence
 \[ 0 \longrightarrow R(A/F) \longrightarrow \Sel_{p^{\infty}}(A/F)
 \longrightarrow \bigoplus_{v|p}A(F_v)\otimes\Qp/\Zp. \]
 In particular, the definition of the fine
Selmer group does not depend on the choice of the set $S$. \el

\bpf Let $S$ be a finite set of primes of $F$ which contains the
primes above $p$, the primes of bad reduction of $A$ and the
archimedean primes. Then by \cite[Chap. I, Corollary 6.6]{Mi}, we
have the following description of the Selmer group
 \[ \Sel_{p^{\infty}}(A/F) = \ker\Big(H^1(G_S(F),A[p^{\infty}])\longrightarrow
\bigoplus_{v\in S}H^1(F_v, A)[p^{\infty}]\Big). \]
 Combining this description with the definition of the fine Selmer
group and an easy diagram-chasing argument, we obtain the required
exact sequence (noting that $A(F_v)\otimes\Qp/\Zp =0$ for $v\nmid
p$). \epf

\br In \cite{Wu}, Wuthrich used the exact sequence in the lemma for
the definition of the fine Selmer group.
 \er

We end the section with the following simple lemma which gives a
lower bound for the $p$-rank of the fine Selmer group in terms of the
$p$-rank of the $S$-class group. This will be used in Sections
\ref{unboundness} and \ref{unramified pro-p}.

\bl \label{lower bound}
  Let $A$ be an abelian variety defined over a
number field $F$. Suppose that $A(F)[p]\neq 0$. Then we have
 \[r_p\big(R(A/F)\big) \geq r_p(\Cl_S(F))r_p(A(F)[p])-2d, \]
 where $d$ denotes the dimension of the abelian variety $A$.
\el

\bpf
 Let $H_S$ be the $p$-Hilbert $S$-class field of $F$ which, by definition, is the maximal abelian unramified $p$-extension
 of $F$ in which all primes in $S$ split completely. Consider the
 following diagram
 \[ \entrymodifiers={!! <0pt, .8ex>+} \SelectTips{eu}{}
\xymatrix{
  0 \ar[r] & R(A/F) \ar[d]^{\al} \ar[r] & H^1(G_S(F),A[p^{\infty}]) \ar[d]^{\beta}
  \ar[r]^{} & \displaystyle\bigoplus_{v\in S}H^1(F_v,
A[p^{\infty}]) \ar[d]^{\gamma} \\
  0 \ar[r] & R(A/H_S) \ar[r] & H^1(G_S(H_S),A[p^{\infty}])
  \ar[r] & \displaystyle\bigoplus_{v\in S}\bigoplus_{w|v}H^1(H_{S,w},
A[p^{\infty}])   }
\]
 with exact rows. Here the vertical maps are given by the restriction maps.
Write $\gamma = \oplus_v \gamma_v$, where
 \[ \gamma_v : H^1(F_v,
A[p^{\infty}]) \longrightarrow \bigoplus_{w|v}H^1(H_{S,w},
A[p^{\infty}]). \]
 It follows from the inflation-restriction sequence that $\ker
 \gamma_v =  H^1(G_v,
A(H_{S,v})[p^{\infty}])$, where $G_v$ is the decomposition group of
$\Gal(H_S/F)$ at $v$. On the other hand, by the definition of $H_S$,
all the primes of $F$ in $S$ split completely in $H_S$, and
therefore, we have $G_v=1$ which in turn implies that $\ker \gamma
=0$. Similarly, the inflation-restriction sequence gives the
equality $\ker \beta = H^1(\Gal(H_S/F), A(H_S)[p^{\infty}])$.
Therefore, we obtain an injection
 \[
 H^1(\Gal(H_S/F), A(H_S)[p^{\infty}])\hookrightarrow R(A/F) \]
It follows from this injection that we have
\[ r_p(R(A/F)) \geq r_p\big(H^1(\Gal(H_S/F), A(H_S)[p^{\infty}])\big). \]
By Lemma \ref{cohomology rank inequalities}, the latter quantity is
greater or equal to
\[ h_1(\Gal(H_S/F))r_p(A(F)[p^{\infty}]) - 2d. \]
By class field theory, we have $\Gal(H_S/F)\cong \Cl_S(F)$, and
therefore,
\[ h_1(\Gal(H_S/F)) = r_p(\Cl_S(F)/p) = r_p(\Cl_S(F)),\]
where the last equality follows from the fact that $\Cl_S(F)$ is
finite. The required estimate is now established (and noting that
$r_p(A(F)[p]) = r_p(A(F)[p^{\infty}]))$.
 \epf

\br
 Since the fine Selmer group is contained in the classical Selmer group
(cf. Lemma \ref{indep of S}), the above estimate also gives a lower
bound for the classical Selmer group. \er

\section{Growth of fine Selmer groups in a $\Zp$-extension}
\label{cyclotomic Zp-extension}

As before, $p$ denotes an odd prime. In this section, $F_{\infty}$
will always denote a fixed $\Zp$-extension of $F$. We will denote
$F_n$ to be the subfield of $F_{\infty}$ such that $[F_n : F] =
p^n$. If $S$ is a finite set of primes of $F$, we denote by $S_f$
the set of finite primes in $S$.
\medskip

We now state the main theorem of this section which compares the
growth of the fine Selmer groups and the growth of the class groups
in the $\Zp$-extension of $F$. To simplify our discussion, we will
assume that $A[p]\subseteq A(F)$.
\medskip

\bt \label{asymptotic compare}
 Let $A$ be a $d$-dimensional abelian variety defined over a number field
 $F$. Let $F_{\infty}$ be a fixed $\Zp$-extension of $F$ such that
the primes of $F$ above $p$ and the bad reduction primes of $A$
decompose finitely in $F_{\infty}/F$. Furthermore, we assume that
$A[p]\subseteq A(F)$.  Then we have
 \[ \Big|r_p(R(A/F_n)) - 2dr_p(\Cl(F_n))\Big|=
 O(1).\]
 \et

In preparation for the proof of the theorem, we require a few
lemmas.

\bl
  Let $F_{\infty}$ be a $\Zp$-extension of $F$ and
let $F_n$ be the subfield of $F_{\infty}$ such that $[F_n : F] =
p^n$. Let $S$ be a given finite set of primes of $F$ which contains
all the primes above $p$ and the archimedean primes. Suppose that
all the primes in $S_f$ decompose finitely in $F_{\infty}/F$. Then
we have
 \[ \Big|r_p(\Cl(F_n)) - r_p(\Cl_S(F_n))\Big|=
  O(1). \]
\el

\bpf
 For each $F_n$, we write $S_f(F_n)$ for the set of finite primes of
$F_n$ above $S_f$. For each $n$, we have the following exact
sequence (cf. \cite[Lemma 10.3.12]{NSW})
 \[ \Z^{|S_f(F_n)|} \longrightarrow \Cl(F_n)\longrightarrow \Cl_S(F_n)
 \longrightarrow 0. \]
Denote by $C_n$  the kernel of $\Cl(F_n)\longrightarrow
\Cl_S(F_n)$. Note that $C_n$ is finite, since it is contained in
$\Cl(F_n)$. Also, it is clear from the above exact sequence that
$r_p(C_n) \leq |S_f(F_n)|$ and $r_p(C_n/p) \leq |S_f(F_n)|$. By
Lemma \ref{estimate lemma}, we have
\[ \Big| r_p(\Cl(F_n)) - r_p(\Cl_S(F_n)) \Big| \leq  3|S_f(F_n)| =
 O(1), \]
 where the last equality follows from the assumption that
all the primes in $S_f$ decompose finitely in $F_{\infty}/F$. \epf
\medskip

Before stating the next lemma, we introduce the $p$-fine Selmer
group of an abelian variety $A$. Let $S$ be a finite set of primes
of $F$ which contains the primes above $p$, the primes of bad
reduction of $A$ and the archimedean primes. Then the $p$-fine
Selmer group (with respect to $S$) is defined to be
\[ R_S(A[p]/F) = \ker\Big(H^1(G_S(F),A[p])\longrightarrow \bigoplus_{v\in S}H^1(F_v,
A[p])\Big). \] Note that the $p$-fine Selmer group may be dependent
on $S$. In fact, as we will see in the proof of Theorem
\ref{asymptotic compare} below, when $F = F(A[p])$, we have
$R_S(A[p]/F) = \Cl_S(F)[p]^{2d}$, where the latter group is clearly
dependent on $S$.

We can now state the following lemma which compare the growth of
$r_p(R_S(A[p]/F_n))$ and $r_p(R(A/F_n))$.

\bl
  Let $F_{\infty}$ be a $\Zp$-extension of $F$ and
let $F_n$  be the subfield of $F_{\infty}$ such that $[F_n : F] =
p^n$. Let $A$ be an abelian variety defined over $F$. Let $S$ be a
finite set of primes of $F$ which contains the primes above $p$, the
primes of bad reduction of $A$ and the archimedean primes.  Suppose
that all the primes in $S_f$ decompose finitely in $F_{\infty}/F$.
Then we have
 \[ \Big|r_p(R_S(A[p]/F_n)) -r_p(R(A/F_n))\Big| =
  O(1).\]
\el

\bpf
 We have a commutative diagram
 \[ \entrymodifiers={!! <0pt, .8ex>+} \SelectTips{eu}{}
\xymatrix{
  0 \ar[r] & R_S(A[p]/F_n) \ar[d]^{s_n} \ar[r] & H^1(G_S(F_n),A[p]) \ar[d]^{h_n}
  \ar[r]^{} & \displaystyle\bigoplus_{v_n\in S(F_n)}H^1(F_{n,v_n},
A[p]) \ar[d]^{g_n} \\
  0 \ar[r] & R(A/F_n)[p] \ar[r] & H^1(G_S(F_n),A[p^{\infty}])[p]
  \ar[r] & \displaystyle\bigoplus_{v_n\in S(F_n)} H^1(F_{n,v_n},
A[p^{\infty}])[p]   }
\]
 with exact rows. It is an easy exercise to show that the maps
$h_n$ and $g_n$ are surjective, that $\ker h_n =
A(F_n)[p^{\infty}]/p$ and that
 \[ \ker g_n = \displaystyle\bigoplus_{v_n\in S(F_n)}A(F_{n,v_n})[p^{\infty}]/p.
 \]
 Since we are assuming $p$ is odd, we have $r_p(\ker g_n)\leq 2d|S_f(F_n)|$.
By an application of Lemma \ref{estimate lemma}, we have
\[ \ba{rl} \Big|  r_p(R_S(A[p]/F_n))) - r_p(R(A/F_n)) \Big|\!
&\leq ~2r_p(\ker s_n) + r_p(\coker s_n) \\
& \leq ~ 2r_p(\ker h_n) + r_p(\ker g_n) \\
& \leq ~ 4d + 2d|S_f(F_n)| = O(1),\\
  \ea \] where the last equality
follows from the assumption that all the primes in $S_f$ decompose
finitely in $F_{\infty}/F$. \epf

We are in the position to prove our theorem.

\bpf[Proof of Theorem \ref{asymptotic compare}]
  Let $S$ be the
finite set of primes of $F$ consisting precisely of the primes above
$p$, the primes of bad reduction of $A$ and the archimedean primes.
By the hypothesis $A[p]\subseteq A(F)$ ($\subseteq A(F_n)$) of the
theorem, we have $A[p] \cong (\Z/p)^{2d}$ as $G_S(F_n)$-modules.
Therefore, we have $H^1(G_S(F_n),A[p]) = \Hom(G_S(F_n),A[p])$. We
have similar identification for the local cohomology groups, and it
follows that
$$
R_S(A[p]/F_n) = \Hom(\Cl_S(F_n),A[p])\cong \Cl_S(F_n)[p]^{2d}
$$
as abelian groups. Hence we have $r_p(R_S(A[p]/F_n)) = 2d
r_p(\Cl_S(F_n))$. The conclusion of the theorem is now immediate
from this equality and the above two lemmas. \epf

\bc \label{asymptotic compare corollary}
 Retain the notations and assumptions of Theorem \ref{asymptotic compare}.
  Then we have
 \[ r_p(R(A/F_n)) = O(1)\]
 if and only if
 \[ r_p(\Cl(F_n)) = O(1).\]
 \ec

For the remainder of the section, $F_{\infty}$ will be taken to be
the cyclotomic $\Zp$-extension of $F$. As before, we denote by $F_n$
the subfield of $F_{\infty}$ such that $[F_n : F] = p^n$. Denote by
$X_{\infty}$ the Galois group of the maximal abelian unramified
pro-$p$ extension of $F_{\infty}$ over $F_{\infty}$. A well-known
conjecture of Iwasawa asserts that $X_{\infty}$ is finitely
generated over $\Zp$ (see \cite{Iw, Iw2}). We will call this
conjecture the \textit{Iwasawa $\mu$-invariant conjecture} for
$F_{\infty}$. By \cite[Proposition 13.23]{Wa}, this is also
equivalent to saying that $r_p(\Cl(F_n)/p)$ is bounded independently
of $n$. Now, by the finiteness of class groups, we have
$r_p(\Cl(F_n))= r_p(\Cl(F_n)/p)$. Hence the Iwasawa $\mu$-invariant
conjecture is equivalent to saying that $r_p(\Cl(F_n))$ is bounded
independently of $n$.
\medskip

We consider the analogous situation for the fine Selmer group. Define
$R(A/F_{\infty}) = \ilim_nR(A/F_n)$ and denote by $Y(A/F_{\infty})$
 the Pontryagin dual of $R(A/F_{\infty})$. We may now recall the
following conjecture which was first introduced in \cite{CS}.

\bigskip \noindent \textbf{Conjecture A.} For any number field $F$,
$Y(A/F_{\infty})$ is a finitely generated $\Zp$-module, where
$F_{\infty}$ is the cyclotomic $\Zp$-extension of $F$.

\medskip
 We can now give the proof of \cite[Theorem 3.4]{CS}. For another
 alternative approach, see \cite{JhS, LimFine}.

\bt \label{Coates-Sujatha}
 Let $A$ be a $d$-dimensional abelian variety defined over a number field
 $F$ and let $F_{\infty}$ be the cyclotomic $\Zp$-extension of $F$.
 Suppose that $F(A[p])$ is a finite $p$-extension of $F$.
  Then Conjecture A holds for $A$ over $F_{\infty}$ if and only if
  the Iwasawa $\mu$-invariant conjecture holds for $F_{\infty}$.
 \et

\bpf
 Now if $L'/L$ is a finite $p$-extension, it follows from \cite[Theorem 3]{Iw}
that the Iwasawa $\mu$-invariant conjecture holds for $L_{\infty}$
if and only if the Iwasawa $\mu$-invariant conjecture holds for
$L'_{\infty}$. On the other hand, it is not difficult to show that
the map
\[ Y(A/L'_{\infty})_{G}\longrightarrow Y(A/L_{\infty})\]
has finite kernel and cokernel, where $G=\Gal(L'/L)$. It follows
from this observation that Conjecture A holds for $A$ over
$L_{\infty}$ if and only if $Y(A/L'_{\infty})_{G}$ is finitely
generated over $\Zp$. Since $G$ is a $p$-group, $\Zp[G]$ is local
with a unique maximal (two-sided) ideal $p\Zp[G]+I_G$, where $I_G$
is the augmentation ideal (see \cite[Proposition 5.2.16(iii)]{NSW}).
It is easy to see from this that
 \[ Y(A/L'_{\infty})/\m \cong Y(A/L'_{\infty})_{G}/pY(A/L'_{\infty})_{G}.
  \] Therefore, Nakayama's lemma
for $\Zp$-modules tells us that $Y(A/L'_{\infty})_{G}$ is finitely
generated over $\Zp$ if and only if $Y(A/L'_{\infty})/\m$ is finite.
On the other hand, Nakayama's lemma for $\Zp[G]$-modules tells us
that $Y(A/L'_{\infty})/\m$ is finite if and only if
$Y(A/L'_{\infty})$ is finitely generated over $\Zp[G]$. But since
$G$ is finite, the latter is equivalent to $Y(A/L'_{\infty})$ being
finitely generated over $\Zp$. Hence we have shown that Conjecture A
holds for $A$ over $L_{\infty}$ if and only if Conjecture A holds
for $A$ over $L_{\infty}'$.
\medskip

Therefore, replacing $F$ by $F(A[p])$, we may assume that
$A[p]\subseteq A(F)$. Write $\Ga_n = \Gal(F_{\infty}/F_n)$. Consider
the following commutative diagram
 \[ \entrymodifiers={!! <0pt, .8ex>+} \SelectTips{eu}{}
\xymatrix{
  0 \ar[r] & R(A/F_n) \ar[d]^{r_n} \ar[r] &
  H^1(G_S(F_n),A[p^{\infty}]) \ar[d]^{f_n}
  \ar[r]^{} & \displaystyle\bigoplus_{v_n\in S(F_n)}H^1(F_{n,v_n},
A[p^{\infty}]) \ar[d]^{\gamma_n} \\
  0 \ar[r] & R(A/F_{\infty})^{\Ga_n} \ar[r] &
  H^1(G_S(F_{\infty}),A[p^{\infty}])^{\Ga_n}
  \ar[r] & \Big(\displaystyle \ilim_n\bigoplus_{v_n\in S(F_n)}
  H^1(F_{n,v_n}, A[p^{\infty}])\Big)^{\Ga_n}   }
\]
 with exact rows, and the vertical maps  given by the restriction maps.
It is an easy exercise to show that $r_p(\ker f_n) \leq 2d$,
$r_p(\ker \gamma_n) \leq 2d|S_f(F_n)|$, and that $f_n$ and
$\gamma_n$ are surjective. It then follows from these estimates and
Lemma \ref{estimate lemma} that we have
\[ \Big|  r_p\big(R(A/F_n))\big) -
r_p\big(R(A/F_{\infty})^{\Ga_n}\big) \Big| = O(1).
\] Combining this observation with \cite[Lemma 13.20]{Wa}, we have
that Conjecture A holds for $A$ over $F_{\infty}$ if and only if
$r_p(R(A/F_n))=O(1)$. The conclusion of the theorem is now immediate
from Corollary \ref{asymptotic compare corollary}. \epf

\section{Unboundedness of fine Selmer groups in $\Z/p$-extensions} \label{unboundness}

In this section, we will study the question of unboundedness of fine
Selmer groups in $\Z/p$-extensions. We first recall the case of
class groups. Since for a number field $L$, the $S$-class
group $\Cl_{S}(L)$ is finite, we have
$r_p(\Cl_{S}(L)) = \dim_{\Z/p}\big(\Cl_S(L)/p\big)$.

\bp \label{class Z/p}
 Let $S$ be a finite set of primes of $F$ which contains all the
the archimedean primes. Then there exists a sequence $\{L_n\}$ of
distinct number fields such that each $L_n$ is a $\Z/p$-extension of
$F$ and such that
\[ r_p(\Cl_{S}(L_n)) \geq n \] for every $n \geq 1$.
\ep

\bpf
 Denote $r_1$ (resp. $r_2$) be the number of real places of $F$
(resp. the number of the pairs of complex places of $F$). Let $S_1$
be a set of primes of $F$ which contains $S$ and such that
\[ |S_1| \geq |S| + r_1 + r_2 + \delta+1. \]
Here $\delta = 1$ if $F$ contains a primitive $p$-root of unity, and
0 otherwise. By the theorem of Grunwald-Wang (cf. \cite[Theorem
9.2.8]{NSW}), there exists a $\Z/p$-extension $L_1$ of $F$ such that
$L_1/F$ is ramified at all the finite primes of $S_1$ and unramified
outside $S_1$. By \cite[Proposition 10.10.3]{NSW}, we have
\[ r_p(\Cl_{S}(L_1)) \geq |S_1| - |S| -r_1 - r_2  -\delta
\geq 1.\]
  Choose $S_2$ to be a set of primes of $F$ which contains
$S_1$ (and hence $S_0$) and which has the property that
\[ |S_2| \geq |S_1| + 1 \geq |S| + r_1 + r_2 + \delta+2. \]
 By the theorem of Grunwald-Wang, there exists a $\Z/p$-extension
$L_2$ of $F$ such that $L_2/F$ is ramified at all the finite primes
of $S_2$ and unramified outside $S_2$. In particular, the fields
$L_1$ and $L_2$ are distinct. By an application of \cite[Proposition
10.10.3]{NSW} again, we have
\[ r_p(\Cl_{S}(L_2)) \geq |S_2| - |S| -r_1 - r_2  -\delta
\geq 2.\]
 Note that since there are infinitely many primes in $F$, we can always continue
the above process iteratively. Also, it is clear from our choice of
$L_n$, they are mutually distinct. Therefore, we have the required
conclusion. \epf

For completeness and for ease of later comparison, we record the
following folklore result.

\bt Let $F$ be a number field. Then we have
 \[ \sup\{ r_p\big(\Cl(L)\big)~ |~
 \mbox{L/F is a cyclic extension of degree p}\} = \infty \]\et

\bpf
 Since $\Cl(L)$ surjects onto $\Cl_S(L)$, the theorem follows from
 the preceding proposition.
\epf

We now record the analogous statement for the fine Selmer groups.

\bt \label{theorem Z/p} Let $A$ be an abelian variety defined over a
number field $F$. Suppose that $A(F)[p]\neq 0$. Then we have
 \[ \sup\{ r_p\big(R(A/L)\big)~ |~
 \mbox{L/F is a cyclic extension of degree p}\} = \infty \]\et

\bpf
 This follows immediately from combining Lemma \ref{lower bound} and
 Proposition \ref{class Z/p}. \epf

In the case, when $A(F)[p]=0$, we have the following weaker
statement.

\bc \label{theorem Z/p corollary} Let $A$ be a $d$-dimensional
abelian variety defined over a number field $F$. Suppose that
$A(F)[p]=0$. Define
\[m = \min\{[K:F]~|~ A(K)[p]\neq 0\}.\]
Then we have
 \[ \sup\{ r_p\big(R(A/L)\big)~ |~
 \mbox{L/F is an extension of degree pm}\} = \infty \]\ec

\bpf
 This follows from an application of the previous theorem to
 the field $K$.  \epf

\br Clearly $1< m\leq |\mathrm{GL}_{2d}(\Z/p)| = (p^{2d}-1)(p^{2d}
-p)\cdots (p^{2d}-p^{2d-1})$. In fact, we can even do
better\footnote{We thank Christian Wuthrich for pointing this out to
us.}. Write $G=\Gal(F(A[p])/F)$. Note that this is a subgroup of
$\mathrm{GL}_{2d}(\Z/p)$. Let $P$ be a nontrivial point in $A[p]$
and denote by $H$ the subgroup of $G$ which fixes $P$. Set $K =
F(A[p])^H$. It is easy to see that $[K:F] = [G:H] = |O_G(P)|$, where
$O_G(P)$ is the orbit of $P$ under the action of $G$. Since $O_G(P)$
is contained in $A[p]\setminus\{0\}$, we have $m \leq [K:F] =
|O_G(P)| \leq p^{2d}-1$. \er

As mentioned in the introductory section, analogous result to the
above theorem for the classical Selmer groups have been studied (see
\cite{Ba, Br, Ce, K, KS, Mat06, Mat09}). Since the fine Selmer group
is contained in the classical Selmer group (cf. Lemma \ref{indep of
S}), our result recovers the above mentioned results (under our
hypothesis). We note that the work of \cite{Ce} also considered the
cases of a global field of positive characteristic. We should also
mention that in \cite{ClS, Cr}, they have even established the
unboundness of $\Sha(A/L)$ over $\Z/p$-extensions $L$ of $F$ in
certain cases (see also \cite{K, Mat09} for some other related
results in this direction). In view of these results on $\Sha(A/L)$,
one may ask for analogous results for a `fine' Shafarevich-Tate
group.
\medskip

Wuthrich \cite{Wu} introduces such a group as follows. One first
defines a `fine' Mordell-Weil group $M_{p^{\infty}}(A/L)$ by the
exact sequence
$$
 0 \longrightarrow\ M_{p^\infty}(A/L) \longrightarrow \ A(L) \otimes \Q_p/\Z_p
  \longrightarrow\ \displaystyle\bigoplus_{v|p}  A(L_v) \otimes \Q_p/\Z_p.
  $$
Then, the `fine' Shafarevich-Tate group is defined by the exact sequence
$$
0 \ \longrightarrow\ M_{p^\infty}(A/L)\ \longrightarrow\ R_{p^\infty}(A/L)\ \longrightarrow\ \Zhe_{p^\infty}(A/L)\ \longrightarrow
\ 0.
$$
In fact, it is not difficult to show that $\Zhe(A/L)$ is contained
in the ($p$-primary) classical Shafarevich-Tate group (see loc.
cit.).
 One may therefore think of $\Zhe(A/L)$ as
the `Shafarevich-Tate part' of the fine Selmer group.
\medskip

With this definition in hand, one is naturally led to the following question for which we do not have
an answer at present.

\medskip \noindent \textbf{Question.} Retaining the assumptions of
Theorem \ref{theorem Z/p}, do we also have
\[ \sup\{ r_p\big(\Zhe_{p^\infty}(A/L)\big)~ |~
 \mbox{L/F is a cyclic extension of degree p}\} = \infty ?\]
\medskip

\section{Growth of fine Selmer groups in
infinite unramified pro-$p$ extensions} \label{unramified pro-p}

We introduce an interesting class of infinite unramified extensions
of $F$. Let $S$ be a finite set (possibly empty) of primes in $F$.
As before, we denote the $S$-ideal class group of $F$ by $\Cl_S(F)$.
For the remainder of the section, $F_{\infty}$ will denote the
maximal unramified $p$-extension of $F$ in which all primes in $S$
split completely. Write $\Sigma = \Sigma_F= \Gal(F_{\infty}/F)$, and
let $\{ \Sigma_n\}$ be the derived series of $\Sigma$. For each $n$,
the fixed field $F_{n+1}$ corresponding to $\Sigma_{n+1}$ is the
$p$-Hilbert $S$-class field of $F_n$.
\medskip

Denote by $S_{\infty}$  the collection of infinite primes of $F$, and
define $\delta$ to be 0 if $\mu_p\subseteq F$ and 1 otherwise. Let
$r_1(F)$ and $r_2(F)$ denote the number of real and complex places
of $F$ respectively. It is known that if the  inequality
 \[ r_p(\Cl_S(F)) \geq 2+ 2\sqrt{r_1(F)+ r_2(F) + \delta +
|S\setminus S_{\infty}|}\]
holds, then $\Sigma$ is infinite (see
\cite{GS}, and also \cite[Chap.\ X, Theorem 10.10.5]{NSW}). Stark
posed the question on whether $r_p(\Cl_S(F_n))$ tends to infinity in
an infinite $p$-class field tower as $n$ tends to infinity. By class
field theory, we have $r_p(\Cl_S(F_n)) = h_1(\Sigma_n)$. It then
follows from the theorem of Lubotzsky and Mann \cite{LuM} that
Stark's question is equivalent to whether the group $\Sigma$ is
$p$-adic analytic. By the following conjecture of Fontaine-Mazur
\cite{FM}, one does not expect $\Sigma$ to be an analytic group if
it is infinite.
\medskip\par
\medskip \noindent \textbf{Conjecture} (Fontaine-Mazur) \textit{For
any number field $F$,\, the group $\Sigma_F$ has no infinite
$p$-adic analytic quotient.}
\medskip

Without assuming the Fontaine-Mazur Conjecture, we have the
following unconditional (weaker) result, proven by various authors.

\bt \label{p-adic class tower} Let $F$ be a number field. If the
following inequality
  \[ r_p(\Cl_S(F)) \geq 2+ 2\sqrt{r_1(F)+ r_2(F) +
\delta + |S\setminus S_{\infty}|}\] holds, then the group $\Sigma_F$
is not $p$-adic analytic. \et

\bpf When $S$ is the empty set, this theorem has been proved
independently by Boston \cite{B} and Hajir \cite{Ha}. For a general
nonempty $S$, this is proved in \cite[Lemma 2.3]{Ma}. \epf
\medskip\par
Collecting all the information we have, we obtain the following
result which answers an analogue of Stark's question, namely the
growth of the $p$-rank of the fine Selmer groups.

\bt \label{Fine Sel in class tower}
 Let $A$ be an Abelian
variety of dimension $d$ defined over $F$ and let $S$ be a finite
set of primes which contains the primes above $p$, the primes of bad
reduction of $A$ and the archimedean primes. Let $F_{\infty}$ be the
maximal unramified $p$-extension of $F$ in which all primes of the
given set $S$ split completely, and let $F_n$ be defined as above.
Suppose that
\[ r_p(\Cl_S(F)) \geq 2+ 2\sqrt{r_1(F)+ r_2(F) + \delta +
|S\setminus S_{\infty}|}\] holds, and suppose that $A(F)[p]\neq 0$.
Then the $p$-rank of $R(A/F_n)$ is unbounded as $n$ tends to
infinity. \et

\bpf
 By Lemma \ref{lower bound}, we have
 \[ r_p(R(A/F_n)) \geq  r_p(\Cl_S(F_n))r_p(A(F)[p])-2d. \]
 Now by the hypothesis of the theorem, it follows from Theorem \ref{p-adic class
 tower} that $\Sigma_F$ is not $p$-adic analytic. By the theorem of Lubotzsky and Mann
 \cite{LuM}, this in turn implies that $r_p(\Cl_S(F_n))$ is
 unbounded as $n$ tends to infinity. Hence we also have that $r_p(R(A/F_n))$ is unbounded as $n$ tends to
 infinity (note here we also make use of the fact that $r_p(A(F)[p])\neq 0$ which comes from the hypothesis that $A(F)[p]\neq 0$). \epf
\medskip

\br
 (1) The analogue of the above result for the classical Selmer group has been
established in \cite{LM, Ma}. In particular, our result here refines
(and implies) those proved there.

(2) Let $A$ be an abelian variety defined over $F$ with complex
multiplication by $K$, and suppose that $K\subseteq F$. Let
$\mathfrak{p}$ be a prime ideal of $K$ above $p$. Then one can
define a $\mathfrak{p}$-version of the fine Selmer group replacing
$A[p^{\infty}]$ by $A[\mathfrak{p}^{\infty}]$ in the definition of
the fine Selmer group. The above arguments carry over to establish
the fine version of the results in \cite{MO}. \er
\medskip

\begin{ack}
    We like to thank Christian Wuthrich for his interest and
    comments on the paper. We will also like to thank the anonymous
referee for pointing out some mistakes and for giving various
helpful comments that have improved the exposition of the paper.
            \end{ack}

%\footnotesize

\medskip

\end{document}